# Fourth Order Compact Formulation of Navier-Stokes Equations and Driven Cavity Flow at High Reynolds Numbers

E. Erturk* and C. Gökçöl

*Energy Systems Engineering Department, Gebze Institute of Technology, Gebze, Kocaeli 41400, Turkey*

SUMMARY

A new fourth order compact formulation for the steady 2-D incompressible Navier-Stokes equations is presented. The formulation is in the same form of the Navier-Stokes equations such that any numerical method that solve the Navier-Stokes equations can easily be applied to this fourth order compact formulation. In particular in this work the formulation is solved with an efficient numerical method that requires the solution of tridiagonal systems using a fine grid mesh of $601\times601$. Using this formulation, the steady 2-D incompressible flow in a driven cavity is solved up to Reynolds number of 20,000 with fourth order spatial accuracy. Detailed solutions are presented.

## 1. INTRODUCTION

High Order Compact (HOC) formulations are becoming more popular in Computational Fluid Dynamics (CFD) field of study. Compact formulations provide more accurate solutions in a compact stencil.

In finite differences, a standard three point discretization provides second order spatial accuracy and this type of discretization is very widely used. When a high order spatial discretization is desired, ie. fourth order accuracy, then a five point discretization have to be used. However in a five point discretization there is a complexity in handling the points near the boundaries.

High order compact schemes provide fourth order spatial accuracy in a $3\times3$ stencil and this type of compact formulations do not have the complexity near the boundaries that a standard wide (five point) fourth order formulation would have.

*Correspondence to: ercanerturk@gyte.edu.tr

Download figures, tables, data files, fortran codes and etc. from `http://www.cavityflow.com`. All the numerical solutions presented in this study is available to public in this web site.

Contract/grant sponsor: Gebze Institute of Technology; contract/grant number: BAP-2003-A-22

Dennis and Hudson [1], MacKinnon and Johnson [2], Gupta *et al.* [3], Spotz and Carey [4] and Li *et al.* [5] have demonstrated the efficiency of the high order compact schemes on the streamfunction and vorticity formulation of 2-D steady incompressible Navier-Stokes equations.

In the literature, it is possible to find numerous different type of iterative numerical methods for the Navier-Stokes equations. These numerical methods, however, could not be easily used in HOC schemes because of the final form of the HOC formulations used in [1], [2], [3], [4] and [5]. This fact might be counted as a disadvantage of HOC formulations that the coding stage is rather complex due to the resulting stencil used in these studies. It would be very useful if any numerical method for the solution of Navier-Stokes equations described in books and papers could be easily applied to high order compact (HOC) formulations.

In this study, we will present a new fourth order compact formulation. The difference of this formulation with References [1-5] is not in the way that the fourth order compact scheme is obtained. The main difference, however, is in the way that the final form of the equations are written. The main advantage of this formulation is that, any iterative numerical method used for Navier-Stokes equations, can be easily applied to this new HOC formulation, since the final form of the presented HOC formulation is in the same form with the Navier-Stokes equations. Moreover if someone already have a second order accurate ($\mathcal{O}\Delta x^2$) code for the solution of steady 2-D incompressible Navier-Stokes equations, using the presented formulation, they can easily convert their existing code to fourth order accuracy ($\mathcal{O}\Delta x^4$) by just adding some coefficients into their existing code. With this new compact formulation, we have solved the steady 2-D incompressible driven cavity flow at very high Reynolds numbers using a very fine grid mesh to demonstrate the efficiency of this new formulation.

## 2. FOURTH ORDER COMPACT FORMULATION

In non-dimensional form, steady 2-D incompressible Navier-Stokes equations in streamfunction ($\psi$) and vorticity ($\omega$) formulation are given as

$$\frac{\partial^2\psi}{\partial x^2} + \frac{\partial^2\psi}{\partial y^2} = -\omega \tag{1}$$

$$\frac{1}{Re}\frac{\partial^2\omega}{\partial x^2} + \frac{1}{Re}\frac{\partial^2\omega}{\partial y^2} = \frac{\partial\psi}{\partial y}\frac{\partial\omega}{\partial x} - \frac{\partial\psi}{\partial x}\frac{\partial\omega}{\partial y} \tag{2}$$

where $x$ and $y$ are the Cartesian coordinates and $Re$ is the Reynolds number. For first order and second order derivatives the following discretizations are fourth order accurate

$$\frac{\partial\phi}{\partial x} = \phi_x - \frac{\Delta x^2}{6}\frac{\partial^3\phi}{\partial x^3} + \mathcal{O}(\Delta x^4) \tag{3}$$

$$\frac{\partial^2\phi}{\partial x^2} = \phi_{xx} - \frac{\Delta x^2}{12}\frac{\partial^4\phi}{\partial x^4} + \mathcal{O}(\Delta x^4) \tag{4}$$

where $\phi_x$ and $\phi_{xx}$ are standard second order central discretizations such that

$$\phi_x = \frac{\phi_{i+1} - \phi_{i-1}}{2\Delta x} \tag{5}$$

$$\phi_{xx} = \frac{\phi_{i+1} - 2\phi_i + \phi_{i-1}}{\Delta x^2} \tag{6}$$

If we apply the discretizations in Equations (3) and (4) to Equations (1) and (2), we obtain the following equations

$$\psi_{xx} + \psi_{yy} - \frac{\Delta x^2}{12}\frac{\partial^4 \psi}{\partial x^4} - \frac{\Delta y^2}{12}\frac{\partial^4 \psi}{\partial y^4} + \mathcal{O}(\Delta x^4, \Delta y^4) = -\omega \tag{7}$$

$$\frac{1}{Re}\omega_{xx} + \frac{1}{Re}\omega_{yy} - \frac{1}{Re}\frac{\Delta x^2}{12}\frac{\partial^4 \omega}{\partial x^4} - \frac{1}{Re}\frac{\Delta y^2}{12}\frac{\partial^4 \omega}{\partial y^4} + \mathcal{O}(\Delta x^4, \Delta y^4) = \psi_y\omega_x - \psi_x\omega_y$$

$$-\frac{\Delta y^2}{6}\omega_x\frac{\partial^3 \psi}{\partial y^3} - \frac{\Delta x^2}{6}\psi_y\frac{\partial^3 \omega}{\partial x^3} + \frac{\Delta x^2}{6}\omega_y\frac{\partial^3 \psi}{\partial x^3} + \frac{\Delta y^2}{6}\psi_x\frac{\partial^3 \omega}{\partial y^3} + \mathcal{O}(\Delta x^4, \Delta x^2\Delta y^2, \Delta y^4) \tag{8}$$

In these equations we have third and fourth derivatives ($\partial^3/\partial x^3$, $\partial^3/\partial y^3$, $\partial^4/\partial x^4$ and $\partial^4/\partial y^4$) of streamfunction and vorticity ($\psi$ and $\omega$) variables. In order to find an expression for these derivatives we use Equations (1) and (2). For example, when we take the first and second $x$-derivative ($\partial/\partial x$ and $\partial^2/\partial x^2$) of the streamfunction equation (1) we obtain the following equations

$$\frac{\partial^3 \psi}{\partial x^3} = -\frac{\partial \omega}{\partial x} - \frac{\partial^3 \psi}{\partial x \partial y^2} \tag{9}$$

$$\frac{\partial^4 \psi}{\partial x^4} = -\frac{\partial^2 \omega}{\partial x^2} - \frac{\partial^4 \psi}{\partial x^2 \partial y^2} \tag{10}$$

And also, by taking the first and second $y$-derivative ($\partial/\partial y$ and $\partial^2/\partial y^2$) of the streamfunction equation (1) we obtain the following equations

$$\frac{\partial^3 \psi}{\partial y^3} = -\frac{\partial \omega}{\partial y} - \frac{\partial^3 \psi}{\partial x^2 \partial y} \tag{11}$$

$$\frac{\partial^4 \psi}{\partial y^4} = -\frac{\partial^2 \omega}{\partial y^2} - \frac{\partial^4 \psi}{\partial x^2 \partial y^2} \tag{12}$$

Using standard second order central discretizations given in Table I, these equations ((9), (10), (11) and (12)) can be written as the followings

$$\frac{\partial^3\psi}{\partial x^3} = -\omega_x - \psi_{xyy} + \mathcal{O}(\Delta x^2, \Delta y^2) \quad (13)$$

$$\frac{\partial^4\psi}{\partial x^4} = -\omega_{xx} - \psi_{xxyy} + \mathcal{O}(\Delta x^2, \Delta y^2) \quad (14)$$

$$\frac{\partial^3\psi}{\partial y^3} = -\omega_y - \psi_{xxy} + \mathcal{O}(\Delta x^2, \Delta y^2) \quad (15)$$

$$\frac{\partial^4\psi}{\partial y^4} = -\omega_{yy} - \psi_{xxyy} + \mathcal{O}(\Delta x^2, \Delta y^2) \quad (16)$$

When we substitute Equations (14) and (16) into Equation (7) we obtain the following finite difference equation.

$$\psi_{xx} + \psi_{yy} = -\omega - \frac{\Delta x^2}{12}\omega_{xx} - \frac{\Delta y^2}{12}\omega_{yy} - \left(\frac{\Delta x^2}{12} + \frac{\Delta y^2}{12}\right)\psi_{xxyy} + \mathcal{O}(\Delta x^4, \Delta x^2\Delta y^2, \Delta y^4) \quad (17)$$

We note that, the solution of Equation (17) is also a solution to streamfunction equation (1) with fourth order spatial accuracy. Therefore if we numerically solve Equation (17), the solution we obtain will satisfy the streamfunction equation up to fourth order accuracy.

In order to obtain a fourth order approximation for the vorticity equation (2), we follow the same procedure. When we take the first and second derivatives of the vorticity equation (2) with respect to $x$- and $y$-coordinates we obtain the followings

$$\frac{\partial^3\omega}{\partial x^3} = Re\frac{\partial^2\psi}{\partial x\partial y}\frac{\partial\omega}{\partial x} + Re\frac{\partial\psi}{\partial y}\frac{\partial^2\omega}{\partial x^2} - Re\frac{\partial^2\psi}{\partial x^2}\frac{\partial\omega}{\partial y} - Re\frac{\partial\psi}{\partial x}\frac{\partial^2\omega}{\partial x\partial y} - \frac{\partial^3\omega}{\partial x\partial y^2} \quad (18)$$

$$\begin{aligned}\frac{\partial^4\omega}{\partial x^4} &= Re\frac{\partial^3\psi}{\partial x^2\partial y}\frac{\partial\omega}{\partial x} + Re\frac{\partial^2\psi}{\partial x\partial y}\frac{\partial^2\omega}{\partial x^2} + Re\frac{\partial^2\psi}{\partial x\partial y}\frac{\partial^2\omega}{\partial x^2} + Re\frac{\partial\psi}{\partial y}\frac{\partial^3\omega}{\partial x^3} \\ &\quad -Re\frac{\partial^3\psi}{\partial x^3}\frac{\partial\omega}{\partial y} - Re\frac{\partial^2\psi}{\partial x^2}\frac{\partial^2\omega}{\partial x\partial y} - Re\frac{\partial^2\psi}{\partial x^2}\frac{\partial^2\omega}{\partial x\partial y} - Re\frac{\partial\psi}{\partial x}\frac{\partial^3\omega}{\partial x^2\partial y} - \frac{\partial^4\omega}{\partial x^2\partial y^2}\end{aligned} \quad (19)$$

$$\frac{\partial^3\omega}{\partial y^3} = Re\frac{\partial^2\psi}{\partial y^2}\frac{\partial\omega}{\partial x} + Re\frac{\partial\psi}{\partial y}\frac{\partial^2\omega}{\partial x\partial y} - Re\frac{\partial^2\psi}{\partial x\partial y}\frac{\partial\omega}{\partial y} - Re\frac{\partial\psi}{\partial x}\frac{\partial^2\omega}{\partial y^2} - \frac{\partial^3\omega}{\partial x^2\partial y} \quad (20)$$

$$\begin{aligned}\frac{\partial^4\omega}{\partial y^4} &= Re\frac{\partial^3\psi}{\partial y^3}\frac{\partial\omega}{\partial x} + Re\frac{\partial^2\psi}{\partial y^2}\frac{\partial^2\omega}{\partial x\partial y} + Re\frac{\partial^2\psi}{\partial y^2}\frac{\partial^2\omega}{\partial x\partial y} + Re\frac{\partial\psi}{\partial y}\frac{\partial^3\omega}{\partial x\partial y^2} \\ &\quad -Re\frac{\partial^3\psi}{\partial x\partial y^2}\frac{\partial\omega}{\partial y} - Re\frac{\partial^2\psi}{\partial x\partial y}\frac{\partial^2\omega}{\partial y^2} - Re\frac{\partial^2\psi}{\partial x\partial y}\frac{\partial^2\omega}{\partial y^2} - Re\frac{\partial\psi}{\partial x}\frac{\partial^3\omega}{\partial y^3} - \frac{\partial^4\omega}{\partial x^2\partial y^2}\end{aligned} \quad (21)$$

If we substitute Equations (18) and (20) for the third derivatives of vorticity ($\partial^3\omega/\partial x^3$ and $\partial^3\omega/\partial y^3$) into Equations (8), (19) and (21) and also if we substitute Equations (13) and (15) for the third derivatives of streamfunction ($\partial^3\psi/\partial x^3$ and $\partial^3\psi/\partial y^3$) into Equations (8), (19) and (21) and finally if we substitute Equations (19) and (21) for the fourth derivatives of vorticity ($\partial^4\omega/\partial x^4$ and $\partial^4\omega/\partial y^4$) into Equation (8), then we obtain the following equation

$$\omega_{xx} + \omega_{yy} - Re\frac{\Delta x^2}{6}\psi_{xy}\omega_{xx} + Re\frac{\Delta y^2}{6}\psi_{xy}\omega_{yy} + Re^2\frac{\Delta x^2}{12}\psi_y\psi_y\omega_{xx} + Re^2\frac{\Delta y^2}{12}\psi_x\psi_x\omega_{yy} =$$

$$Re\psi_y\omega_x - Re\psi_x\omega_y + Re\left(\frac{\Delta x^2}{12}+\frac{\Delta y^2}{12}\right)\psi_{xxy}\omega_x - Re\left(\frac{\Delta x^2}{12}+\frac{\Delta y^2}{12}\right)\psi_{xyy}\omega_y$$
$$-Re^2\frac{\Delta x^2}{12}\psi_y\psi_{xy}\omega_x + Re^2\frac{\Delta y^2}{12}\psi_x\psi_{yy}\omega_x + Re^2\frac{\Delta x^2}{12}\psi_y\psi_{xx}\omega_y - Re^2\frac{\Delta y^2}{12}\psi_x\psi_{xy}\omega_y$$
$$+Re\left(\frac{\Delta x^2}{12}+\frac{\Delta y^2}{12}\right)\psi_y\omega_{xyy} - Re\left(\frac{\Delta x^2}{12}+\frac{\Delta y^2}{12}\right)\psi_x\omega_{xxy} - Re\frac{\Delta x^2}{6}\psi_{xx}\omega_{xy}$$
$$+Re\frac{\Delta y^2}{6}\psi_{yy}\omega_{xy} + Re^2\left(\frac{\Delta x^2}{12}+\frac{\Delta y^2}{12}\right)\psi_x\psi_y\omega_{xy} - Re\left(\frac{\Delta x^2}{12}-\frac{\Delta y^2}{12}\right)\omega_x\omega_y$$
$$-\left(\frac{\Delta x^2}{12}+\frac{\Delta y^2}{12}\right)\omega_{xxyy} + \mathcal{O}(\Delta x^4, \Delta x^2\Delta y^2, \Delta y^4) \quad (22)$$

Again we note that, the solution of Equation (22) satisfy the vorticity equation (2) with fourth order accuracy.

As the final form of our HOC scheme, we prefer to write Equations (17) and (22) as the following

$$\psi_{xx} + \psi_{yy} = -\omega + A \quad (23)$$

$$\frac{1}{Re}(1+B)\,\omega_{xx} + \frac{1}{Re}(1+C)\,\omega_{yy} = (\psi_y + D)\,\omega_x - (\psi_x + E)\,\omega_y + F \quad (24)$$

where

$$\begin{aligned}
A &= -\frac{\Delta x^2}{12}\omega_{xx} - \frac{\Delta y^2}{12}\omega_{yy} - \left(\frac{\Delta x^2}{12}+\frac{\Delta y^2}{12}\right)\psi_{xxyy}\\
B &= -Re\frac{\Delta x^2}{6}\psi_{xy} + Re^2\frac{\Delta x^2}{12}\psi_y\psi_y\\
C &= Re\frac{\Delta y^2}{6}\psi_{xy} + Re^2\frac{\Delta y^2}{12}\psi_x\psi_x\\
D &= \left(\frac{\Delta x^2}{12}+\frac{\Delta y^2}{12}\right)\psi_{xxy} - Re\frac{\Delta x^2}{12}\psi_y\psi_{xy} + Re\frac{\Delta y^2}{12}\psi_x\psi_{yy}\\
E &= \left(\frac{\Delta x^2}{12}+\frac{\Delta y^2}{12}\right)\psi_{xyy} - Re\frac{\Delta x^2}{12}\psi_y\psi_{xx} + Re\frac{\Delta y^2}{12}\psi_x\psi_{xy}\\
F &= \left(\frac{\Delta x^2}{12}+\frac{\Delta y^2}{12}\right)\psi_y\omega_{xyy} - \left(\frac{\Delta x^2}{12}+\frac{\Delta y^2}{12}\right)\psi_x\omega_{xxy} - \frac{\Delta x^2}{6}\psi_{xx}\omega_{xy}\\
&\quad + \frac{\Delta y^2}{6}\psi_{yy}\omega_{xy} + Re\left(\frac{\Delta x^2}{12}+\frac{\Delta y^2}{12}\right)\psi_x\psi_y\omega_{xy} - \left(\frac{\Delta x^2}{12}-\frac{\Delta y^2}{12}\right)\omega_x\omega_y\\
&\quad - \frac{1}{Re}\left(\frac{\Delta x^2}{12}+\frac{\Delta y^2}{12}\right)\omega_{xxyy}
\end{aligned} \quad (25)$$

We note that, the finite difference Equations (23) and (24) are fourth order accurate $\big(\mathcal{O}(\Delta x^4, \Delta x^2\Delta y^2, \Delta y^4)\big)$ approximation of the streamfunction and vorticity equations (1) and (2). In Equations (23) and (24), however, if $A$, $B$, $C$, $D$, $E$ and $F$ are chosen to be equal to zero then the finite difference Equations (23) and (24) simply become

$$\psi_{xx} + \psi_{yy} = -\omega \tag{26}$$

$$\frac{1}{Re}\omega_{xx} + \frac{1}{Re}\omega_{yy} = \psi_y\omega_x - \psi_x\omega_y \tag{27}$$

Equations (26) and (27) are the standard second order accurate $\big(\mathcal{O}(\Delta x^2, \Delta y^2)\big)$ approximation of the streamfunction and vorticity equations (1) and (2). When we use Equations (23) and (24) for the numerical solution of 2-D steady incompressible Navier-Stokes equations, we can easily switch between second and fourth order accuracy just by using homogeneous values for the coefficients $A$, $B$, $C$, $D$, $E$ and $F$ or by using the expressions defined in Equation (25) in the code.

In Equations (23), (24) and (25) instead of finite difference discretizations, if we substitute for partial derivatives we obtain the following differential equations

$$\frac{\partial^2\psi}{\partial x^2} + \frac{\partial^2\psi}{\partial y^2} = -\omega + A \tag{28}$$

$$\frac{1}{Re}(1+B)\frac{\partial^2\omega}{\partial x^2} + \frac{1}{Re}(1+C)\frac{\partial^2\omega}{\partial y^2} = \left(\frac{\partial\psi}{\partial y} + D\right)\frac{\partial\omega}{\partial x} - \left(\frac{\partial\psi}{\partial x} + E\right)\frac{\partial\omega}{\partial y} + F \tag{29}$$

$$\begin{aligned}
A &= -\frac{\Delta x^2}{12}\frac{\partial^2\omega}{\partial x^2} - \frac{\Delta y^2}{12}\frac{\partial^2\omega}{\partial y^2} - \left(\frac{\Delta x^2}{12} + \frac{\Delta y^2}{12}\right)\frac{\partial^4\psi}{\partial x^2\partial y^2} \\
B &= -Re\frac{\Delta x^2}{6}\frac{\partial^2\psi}{\partial x\partial y} + Re^2\frac{\Delta x^2}{12}\frac{\partial\psi}{\partial y}\frac{\partial\psi}{\partial y} \\
C &= Re\frac{\Delta y^2}{6}\frac{\partial^2\psi}{\partial x\partial y} + Re^2\frac{\Delta y^2}{12}\frac{\partial\psi}{\partial x}\frac{\partial\psi}{\partial x} \\
D &= \left(\frac{\Delta x^2}{12} + \frac{\Delta y^2}{12}\right)\frac{\partial^3\psi}{\partial x^2\partial y} - Re\frac{\Delta x^2}{12}\frac{\partial\psi}{\partial y}\frac{\partial^2\psi}{\partial x\partial y} + Re\frac{\Delta y^2}{12}\frac{\partial\psi}{\partial x}\frac{\partial^2\psi}{\partial y^2} \\
E &= \left(\frac{\Delta x^2}{12} + \frac{\Delta y^2}{12}\right)\frac{\partial^3\psi}{\partial x\partial y^2} - Re\frac{\Delta x^2}{12}\frac{\partial\psi}{\partial y}\frac{\partial^2\psi}{\partial x^2} + Re\frac{\Delta y^2}{12}\frac{\partial\psi}{\partial x}\frac{\partial^2\psi}{\partial x\partial y} \\
F &= \left(\frac{\Delta x^2}{12} + \frac{\Delta y^2}{12}\right)\frac{\partial\psi}{\partial y}\frac{\partial^3\omega}{\partial x\partial y^2} - \left(\frac{\Delta x^2}{12} + \frac{\Delta y^2}{12}\right)\frac{\partial\psi}{\partial x}\frac{\partial^3\omega}{\partial x^2\partial y} - \frac{\Delta x^2}{6}\frac{\partial^2\psi}{\partial x^2}\frac{\partial^2\omega}{\partial x\partial y} \\
&\quad + \frac{\Delta y^2}{6}\frac{\partial^2\psi}{\partial y^2}\frac{\partial^2\omega}{\partial x\partial y} + Re\left(\frac{\Delta x^2}{12} + \frac{\Delta y^2}{12}\right)\frac{\partial\psi}{\partial x}\frac{\partial\psi}{\partial y}\frac{\partial^2\omega}{\partial x\partial y} - \left(\frac{\Delta x^2}{12} - \frac{\Delta y^2}{12}\right)\frac{\partial\omega}{\partial x}\frac{\partial\omega}{\partial y} \\
&\quad - \frac{1}{Re}\left(\frac{\Delta x^2}{12} + \frac{\Delta y^2}{12}\right)\frac{\partial^4\omega}{\partial x^2\partial y^2}
\end{aligned} \tag{30}$$

We note that the numerical solutions of Equations (28) and (29), strictly provided that second order discretizations in Table I are used and also strictly provided that a uniform grid mesh with $\Delta x$ and $\Delta y$ is used, are fourth order accurate to streamfunction and vorticity equations (1) and (2). We prefer to call Equations (28) and (29) Fourth Order Navier-Stokes (FONS) equations. The only difference between FONS equations (28) and (29) and Navier-Stokes (NS) equations (1) and (2) are the coefficients $A$, $B$, $C$, $D$, $E$ and $F$. In fact the NS equations are a subset of the FONS equations. We note that FONS equations (28) and (29) are in the same form with Navier-Stokes (NS) equations (1) and (2), therefore any iterative numerical method (such as SOR, ADI, factorization schemes, pseudo time iterations and etc.) used to solve streamfunction and vorticity equations (1) and (2) can also be easily applied to fourth order equations (28) and (29). Moreover, any existing code that solve the streamfunction and vorticity equations with second order accuracy can easily be modified to provide fourth order accuracy just by adding the coefficients $A$, $B$, $C$, $D$, $E$ and $F$ into the existing code to obtain the solution of FONS equations. Of course, when the coefficients $A$, $B$, $C$, $D$, $E$ and $F$ are added into a second order accurate code to obtain fourth order accuracy, evaluating these coefficients would require extra CPU work. This might be considered as the cost of increasing accuracy from second order to fourth order.

## 3. NUMERICAL METHOD

Recently Erturk *et al.* [6] have presented a new, stable and efficient numerical method that solve the streamfunction and vorticity equations. The numerical method solve the governing steady equations through iterations in the pseudo time. In this study, we will apply the numerical method Erturk *et al.* [6] have proposed, to FONS equations (28) and (29) and solve the steady driven cavity flow with fourth order accuracy. For details about the numerical method, the reader is referred to Erturk *et al.* [6]. When we apply the numerical method to Equations (28) and (29), we obtain the following equations

$$\left(1-\Delta t\frac{\partial^2}{\partial x^2}\right)\left(1-\Delta t\frac{\partial^2}{\partial y^2}\right)\psi^{n+1}=\psi^n+\Delta t\omega^n-\Delta t A^n+\left(\Delta t\frac{\partial^2}{\partial x^2}\right)\left(\Delta t\frac{\partial^2}{\partial y^2}\right)\psi^n \tag{31}$$

$$\begin{aligned}&\left(1-\Delta t(1+B^n)\frac{1}{Re}\frac{\partial^2}{\partial x^2}+\Delta t\left(\frac{\partial\psi}{\partial y}+D\right)^n\frac{\partial}{\partial x}\right)\\ &\qquad\left(1-\Delta t(1+C^n)\frac{1}{Re}\frac{\partial^2}{\partial y^2}-\Delta t\left(\frac{\partial\psi}{\partial x}+E\right)^n\frac{\partial}{\partial y}\right)\omega^{n+1}=\omega^n-\Delta t F^n\\ &\qquad\quad+\left(\Delta t(1+B^n)\frac{1}{Re}\frac{\partial^2}{\partial x^2}-\Delta t\left(\frac{\partial\psi}{\partial y}+D\right)^n\frac{\partial}{\partial x}\right)\\ &\qquad\qquad\left(\Delta t(1+C^n)\frac{1}{Re}\frac{\partial^2}{\partial y^2}+\Delta t\left(\frac{\partial\psi}{\partial x}+E\right)^n\frac{\partial}{\partial y}\right)\omega^n\end{aligned} \tag{32}$$

The solution methodology of these two equations are quite simple. First the streamfunction equation (31) is solved in two steps. For streamfunction equation, a new variable $f$ is defined as the following

$$\left(1-\Delta t\frac{\partial^2}{\partial y^2}\right)\psi^{n+1}=f \tag{33}$$

Using this variable in Equation (31) we obtain the following equation

$$\left(1-\Delta t\frac{\partial^2}{\partial x^2}\right) f = \psi^n + \Delta t\omega^n - \Delta t A^n + \left(\Delta t\frac{\partial^2}{\partial x^2}\right)\left(\Delta t\frac{\partial^2}{\partial y^2}\right)\psi^n \tag{34}$$

In this equation, the only unknown is the variable $f$. We first solve this equation for $f$ by solving a tridiagonal system. After this, when we obtain the value of $f$ at every grid point we solve Equation (33) for streamfunction ($\psi^{n+1}$) by solving another tridiagonal system.

After solving the streamfunction equation (31), we solve the vorticity equation (32). For this, similarly, we introduce a new variable $g$ which is defined as the following

$$\left(1-\Delta t(1+C^n)\frac{1}{Re}\frac{\partial^2}{\partial y^2} - \Delta t\left(\frac{\partial\psi}{\partial x}+E\right)^n\frac{\partial}{\partial y}\right)\omega^{n+1} = g \tag{35}$$

Using this variable in Equation (32), we obtain the following equation

$$\begin{aligned}\left(1-\Delta t(1+B^n)\frac{1}{Re}\frac{\partial^2}{\partial x^2} + \Delta t\left(\frac{\partial\psi}{\partial y}+D\right)^n\frac{\partial}{\partial x}\right) g = \omega^n - \Delta t F^n \\ + \left(\Delta t(1+B^n)\frac{1}{Re}\frac{\partial^2}{\partial x^2} - \Delta t\left(\frac{\partial\psi}{\partial y}+D\right)^n\frac{\partial}{\partial x}\right) \\ \left(\Delta t(1+C^n)\frac{1}{Re}\frac{\partial^2}{\partial y^2} + \Delta t\left(\frac{\partial\psi}{\partial x}+E\right)^n\frac{\partial}{\partial y}\right)\omega^n\end{aligned} \tag{36}$$

In this equation the only unknown is the variable $g$. By solving a tridiagonal system, we obtain the value of $g$ at every grid point. Then we solve Equation (35) for vorticity ($\omega^{n+1}$) by solving another tridiagonal system.

In a compact formulation, the stencil have 3×3 points. The solution at the first diagonal grid points near the corners of the cavity would require the vorticity values at the corner points. However, the corner points are singular points for vorticity. Gupta *et al.* [7] have introduced an explicit asymptotic solution in the neighborhood of sharp corners. Similarly, Störtkuhl *et al.* [8] have presented an analytical asymptotic solutions near the corners of cavity and using finite element bilinear shape functions they also have presented a singularity removed boundary condition for vorticity at the corner points as well as at the wall points. We follow Störtkuhl *et al.* [8] and use the following expression for calculating vorticity values at the wall

$$\frac{1}{3\Delta h^2}\begin{bmatrix} \bullet & \bullet & \bullet \\ \frac{1}{2} & -4 & \frac{1}{2} \\ 1 & 1 & 1 \end{bmatrix}\psi + \frac{1}{9}\begin{bmatrix} \bullet & \bullet & \bullet \\ \frac{1}{2} & 2 & \frac{1}{2} \\ \frac{1}{4} & 1 & \frac{1}{4} \end{bmatrix}\omega = -\frac{V}{h} \tag{37}$$

where $V$ is the speed of the wall which is equal to 1 for the moving top wall and equal to 0 for the three stationary walls. For corner points, we use the following expression for calculating the vorticity values

$$\frac{1}{3\Delta h^2}\begin{bmatrix} \bullet & \bullet & \bullet \\ \bullet & -2 & \frac{1}{2} \\ \bullet & \frac{1}{2} & 1 \end{bmatrix}\psi + \frac{1}{9}\begin{bmatrix} \bullet & \bullet & \bullet \\ \bullet & 1 & \frac{1}{2} \\ \bullet & \frac{1}{2} & \frac{1}{4} \end{bmatrix}\omega = -\frac{V}{2h} \tag{38}$$

where again $V$ is equal to 1 for the upper two corners and it is equal to 0 for the bottom two corners. The reader is referred to Störtkuhl *et al.* [8] for details.

## 4. RESULTS AND DISCUSSIONS

The schematics of the driven cavity flow is given in Figure 1. In this figure the abbreviations BR, BL and TL refer to bottom right, bottom left and top left corners of the cavity, respectively. The number following these abbreviations refer to the vortices that appear in the flow, which are numbered according to size.

For every Reynolds number considered, we have continued our iterations until, in the computational domain both the maximum residual of Equations (23) and (24), which are given as

$$R_\psi = \max\left(\text{abs}\left(\left|\psi_{xx}^{n+1} + \psi_{yy}^{n+1} + \omega^{n+1} - A^{n+1}\right|_{i,j}\right)\right) \tag{39}$$

$$R_\omega = \max\left(\text{abs}\left(\left|\frac{1}{Re}\left(1 + B^{n+1}\right)\omega_{xx}^{n+1} + \frac{1}{Re}\left(1 + C^{n+1}\right)\omega_{yy}^{n+1} - \left(\psi_y^{n+1} + D^{n+1}\right)\omega_x^{n+1} + \left(\psi_x^{n+1} + E^{n+1}\right)\omega_y^{n+1} - F^{n+1}\right|_{i,j}\right)\right) \tag{40}$$

are less than $10^{-10}$. Such a low value is chosen to ensure the accuracy of the solution. At these residual levels, the maximum absolute change in streamfunction value between two time steps, $(\max(|\psi^{n+1} - \psi^n|))$, was in the order of $10^{-16}$ and for vorticity, $(\max(|\omega^{n+1} - \omega^n|))$, it was in the order of $10^{-14}$. Obviously these convergence levels are far more less than satisfactory, however such low values demonstrate the efficiency of the numerical method used in this study which was presented by Erturk *et al.* [6].

Using an efficient numerical method, Erturk *et al.* [6] have clearly shown that numerical solutions of driven cavity flow is computable for $Re$ >10,000 when a grid mesh larger than 256×256 is used. With a grid mesh of 601×601 Erturk *et al.* [6] have solved the cavity flow up to $Re$=21,000 using the numerical method also used in this study. In order to be able to obtain solutions at high Reynolds numbers, following Erturk *et al.* [6], in this study we have used a large grid mesh with 601×601 grids. With this many number of grid points we obtained steady solutions of the cavity flow up to $Re$=20,000 with fourth order accuracy.

Figures 2 to 6 show the streamfunction and vorticity contours of the driven cavity flow between $Re$=1,000 and $Re$=20,000. These figures show the vortices that are formed in the flow field as the Reynolds number increases. From these contour figures, we conclude that the fourth order compact formulation provides very smooth solutions.

In Figure 7 we plot a very enlarged view of the top right corner (where the moving wall moves towards the stationary wall) of the streamfunction contour plot for the highest Reynolds number considered, $Re$=20,000. In this figure the dotted lines show the grid lines. As it is seen in this enlarged figure, fourth order streamfunction contours are very smooth even at the first set of grid points near the corners.

Table II tabulates the streamfunction and vorticity values at the center of the primary and secondary vortices and also the location of the center of these vortices for future references. This table is in good agreement with that of Erturk *et al.* [6].

Using Richardson extrapolation on the solutions obtained with different grid meshes, Erturk *et al.* [6] have presented theoretically fourth and sixth order accurate ($\mathcal{O}\Delta x^4$ and $\mathcal{O}\Delta x^6$) streamfunction and vorticity values at the center of the primary vortex. Table III and IV compares the forth order compact scheme solutions of the streamfunction and the vorticity values at the center of the primary vortex with the fourth order ($\mathcal{O}\Delta x^4$) Richardson extrapolated solutions tabulated in Erturk *et al.* [6]. The present solutions and the solutions of Erturk *et al.* [6] agree with each other.

## 5. CONCLUSIONS

In this study a new fourth order compact formulation is presented. The uniqueness of this formulation is that the final form of the HOC formulation is in the same form of the Navier-Stokes equations such that any numerical method that solve the Navier-Stokes equations can be easily applied to the FONS equations in order to obtain fourth order accurate solutions ($\mathcal{O}\Delta x^4$). Moreover with this formulation, any existing code that solve the Navier-Stokes equations with second order accuracy ($\mathcal{O}\Delta x^2$) can be altered to provide fourth order accurate ($\mathcal{O}\Delta x^4$) solutions just by adding some coefficients into the code at the expense of extra CPU work of evaluating these coefficients.

In this study, the presented fourth order compact formulation is solved with a very efficient numerical method introduced by Erturk *et al.* [6]. Using a fine grid mesh of 601×601, as it was suggested by Erturk *et al.* [6] in order to be able to compute for high Reynolds numbers, the driven cavity flow is solved up to Reynolds numbers of $Re$=20,000. The solutions obtained agree well with previous studies. The presented fourth order accurate compact formulation is proved to be very efficient.

ACKNOWLEDGEMENT

This study was funded by Gebze Institute of Technology with project no BAP-2003-A-22. E. Ertürk is grateful for this financial support.

## REFERENCES

1. S. C. Dennis, J. D. Hudson, Compact $h^4$ Finite Difference Approximations to Operators of Navier-Stokes Type, *Journal of Computational Physics* **85** (1989) 390–416.
2. R. J. MacKinnon, R. W. Johnson, Differential-Equation-Based Representation of Truncation Errors for Accurate Numerical Simulation, *International Journal for Numerical Methods in Fluids* **13** (1991) 739–757.
3. M. M. Gupta, R. P. Manohar, J. W. Stephenson, A Single Cell High Order Scheme for the Convection-Diffusion Equation with Variable Coefficients, *International Journal for Numerical Methods in Fluids* **4** (1984) 641–651.
4. W. F. Spotz, G. F. Carey, High-Order Compact Scheme for the Steady Streamfunction Vorticity Equations, *International Journal for Numerical Methods in Engineering* **38** (1995) 3497–3512.
5. M. Li, T. Tang, B. Fornberg, A Compact Forth-Order Finite Difference Scheme for the Steady Incompressible Navier-Stokes Equations, *International Journal for Numerical Methods in Fluids* **20** (1995) 1137–1151.
6. E. Erturk, T.C. Corke, C. Gokcol, Numerical Solutions of 2-D Steady Incompressible Driven Cavity Flow at High Reynolds Numbers, *International Journal for Numerical Methods in Fluids* **48** (2005) 747–774.
7. M. M. Gupta, R. P. Manohar, B. Noble, Nature of Viscous Flows Near Sharp Corners, *Computers and Fluids* **9** (1981) 379–388.

8. T. Stortkuhl, C. Zenger, S. Zimmer, An Asymptotic Solution for the Singularity at the Angular Point of the Lid Driven Cavity, *International Journal of Numerical Methods for Heat Fluid Flow* **4** (1994) 47–59.

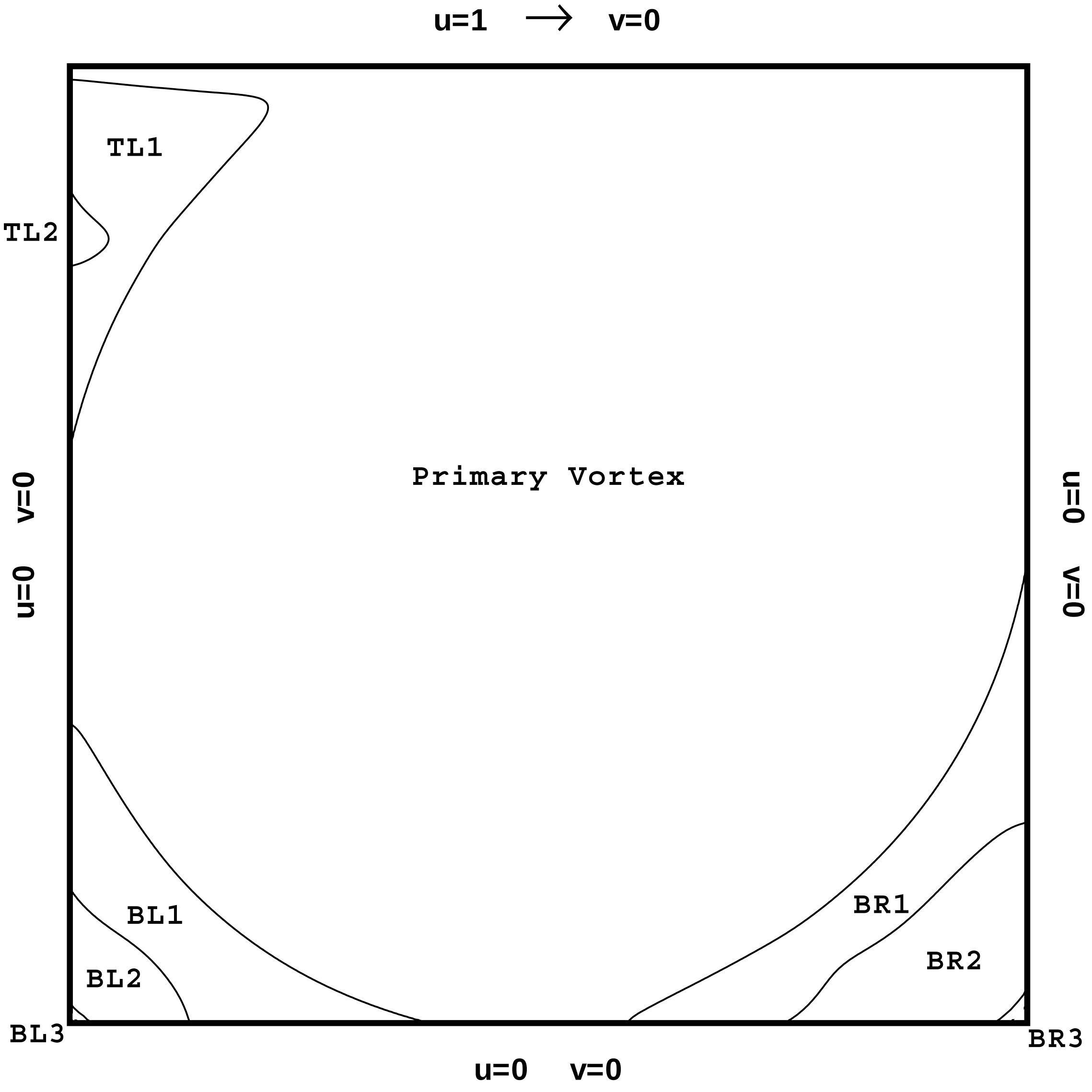

Figure 1. Schematic view of driven cavity flow

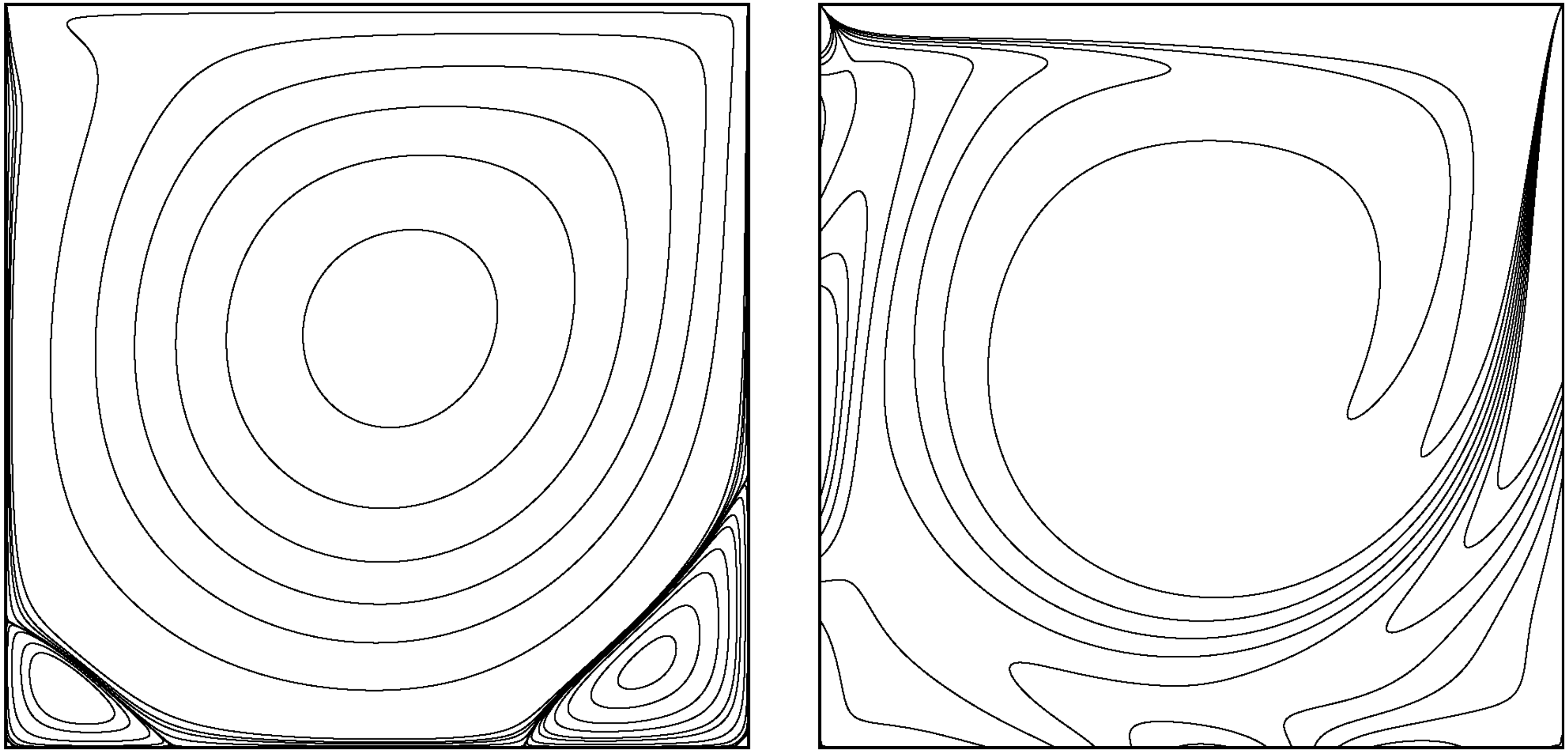

Figure 2.   Streamfunction and vorticity contours for Re=1000

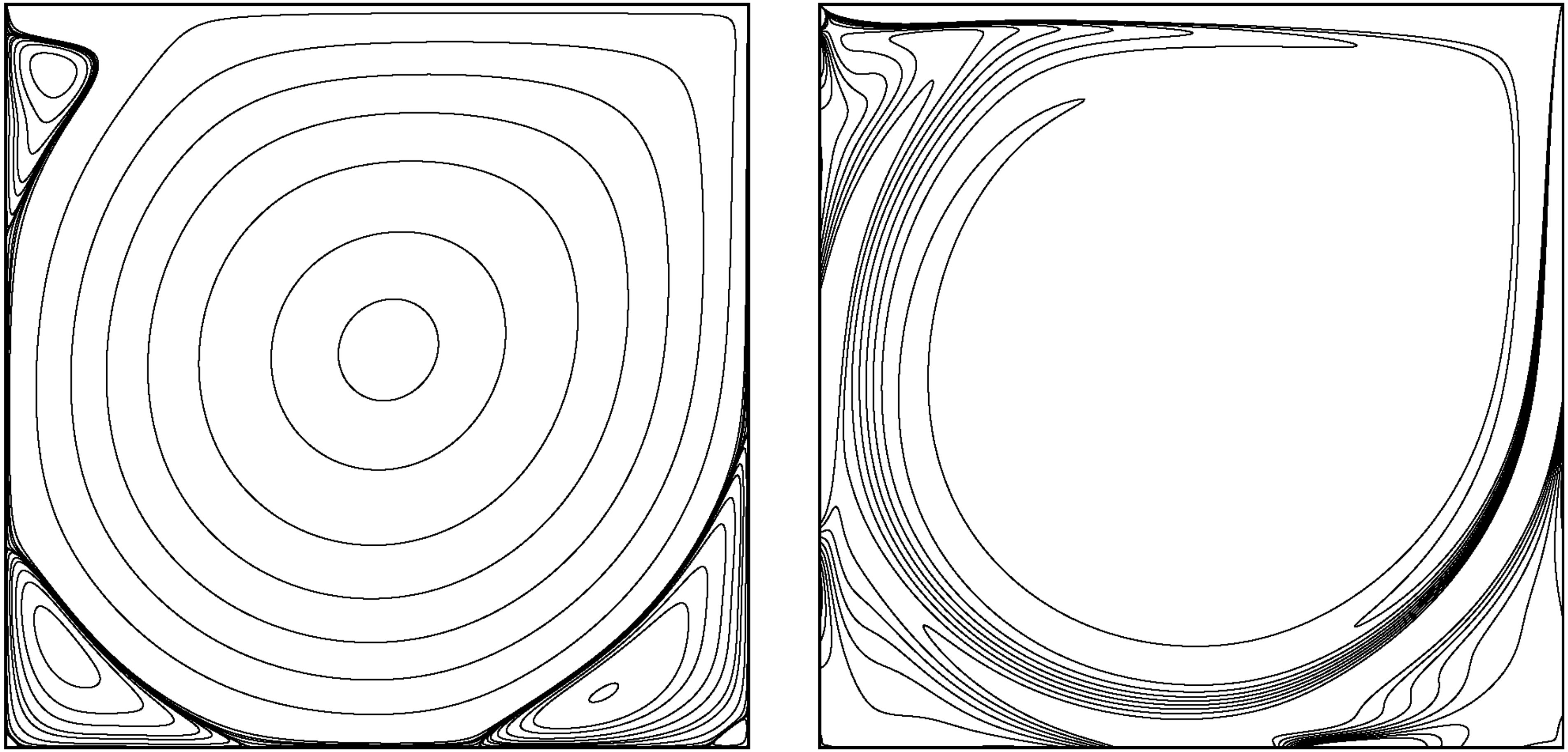

Figure 3. Streamfunction and vorticity contours for Re=5000

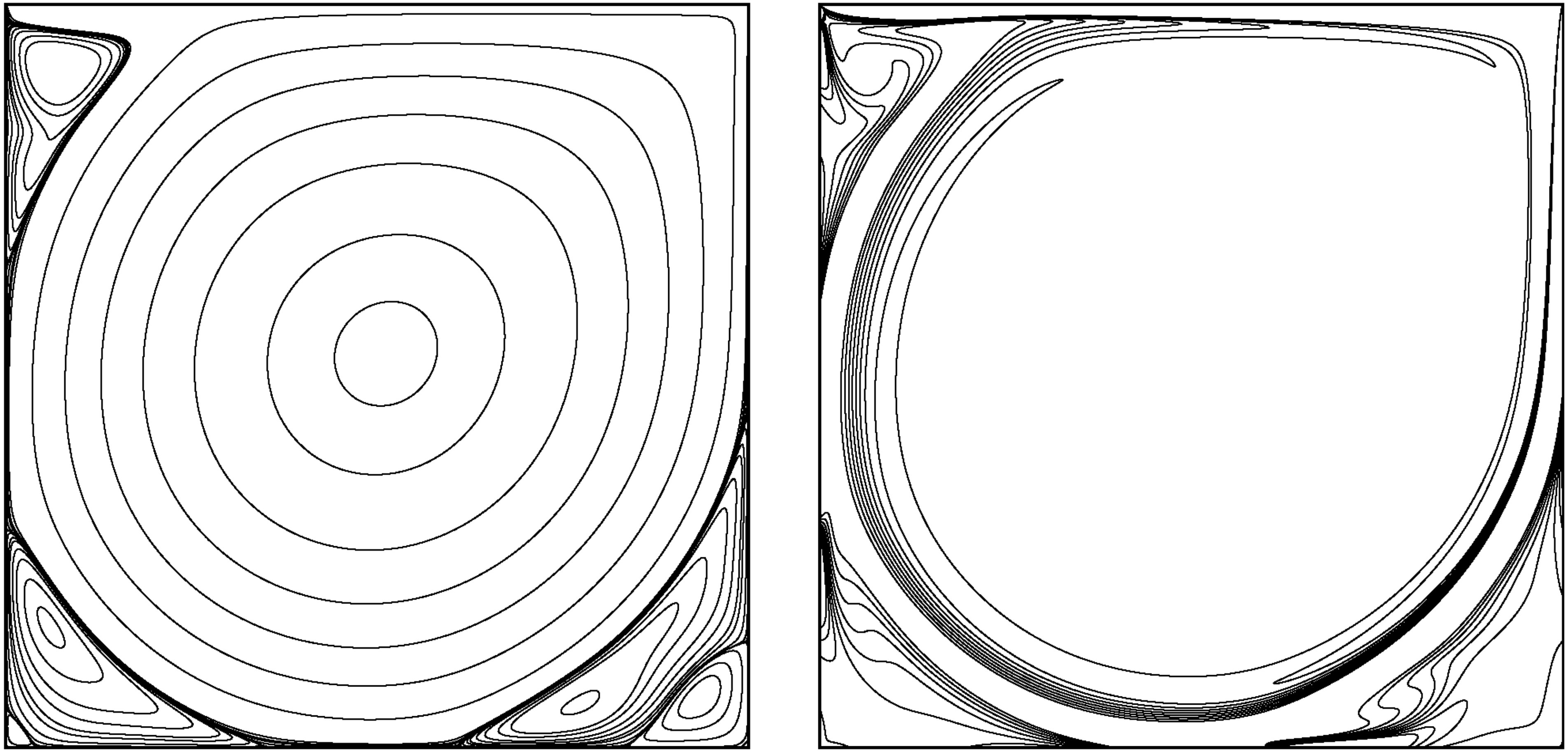

Figure 4. Streamfunction and vorticity contours for Re=10000

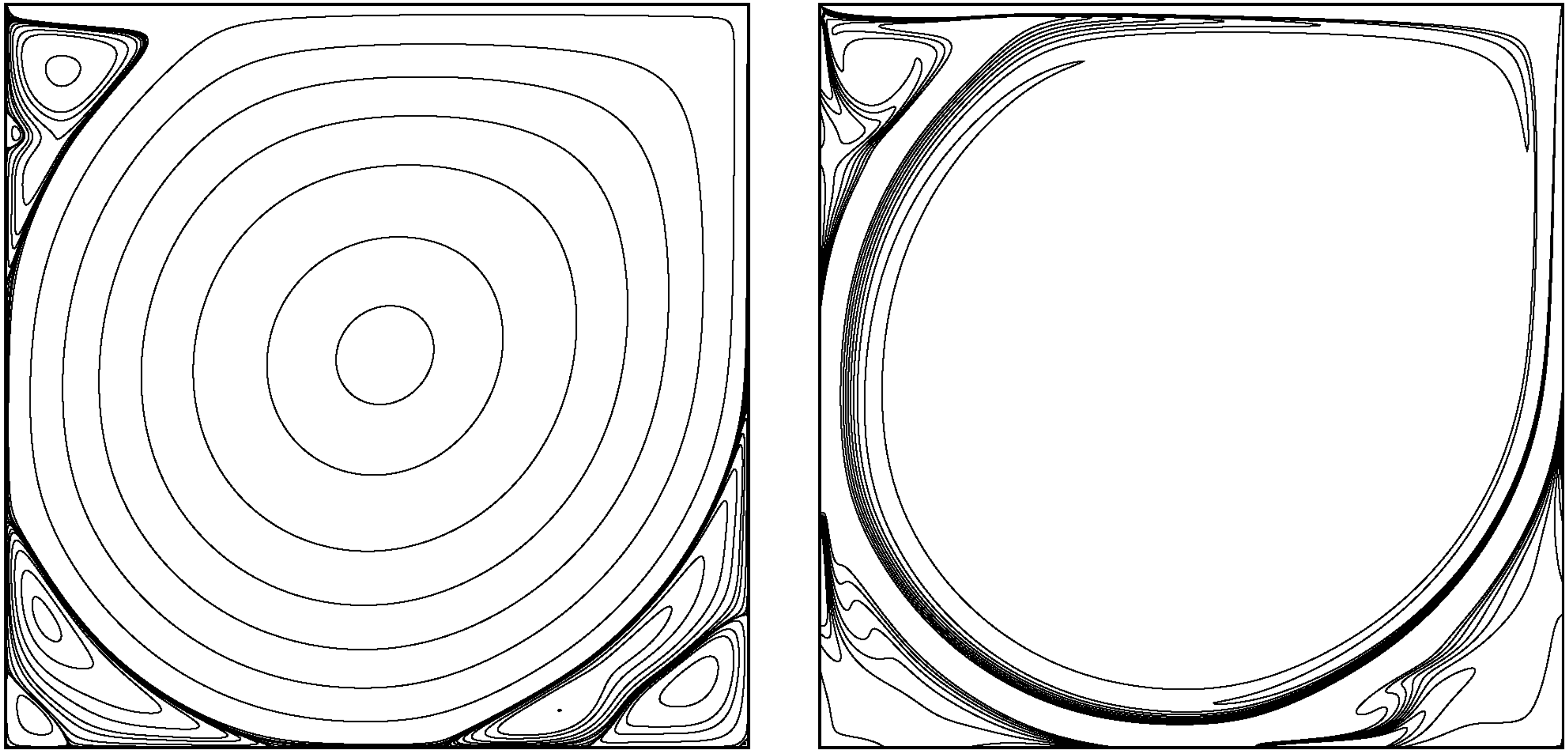

Figure 5. Streamfunction and vorticity contours for Re=15000

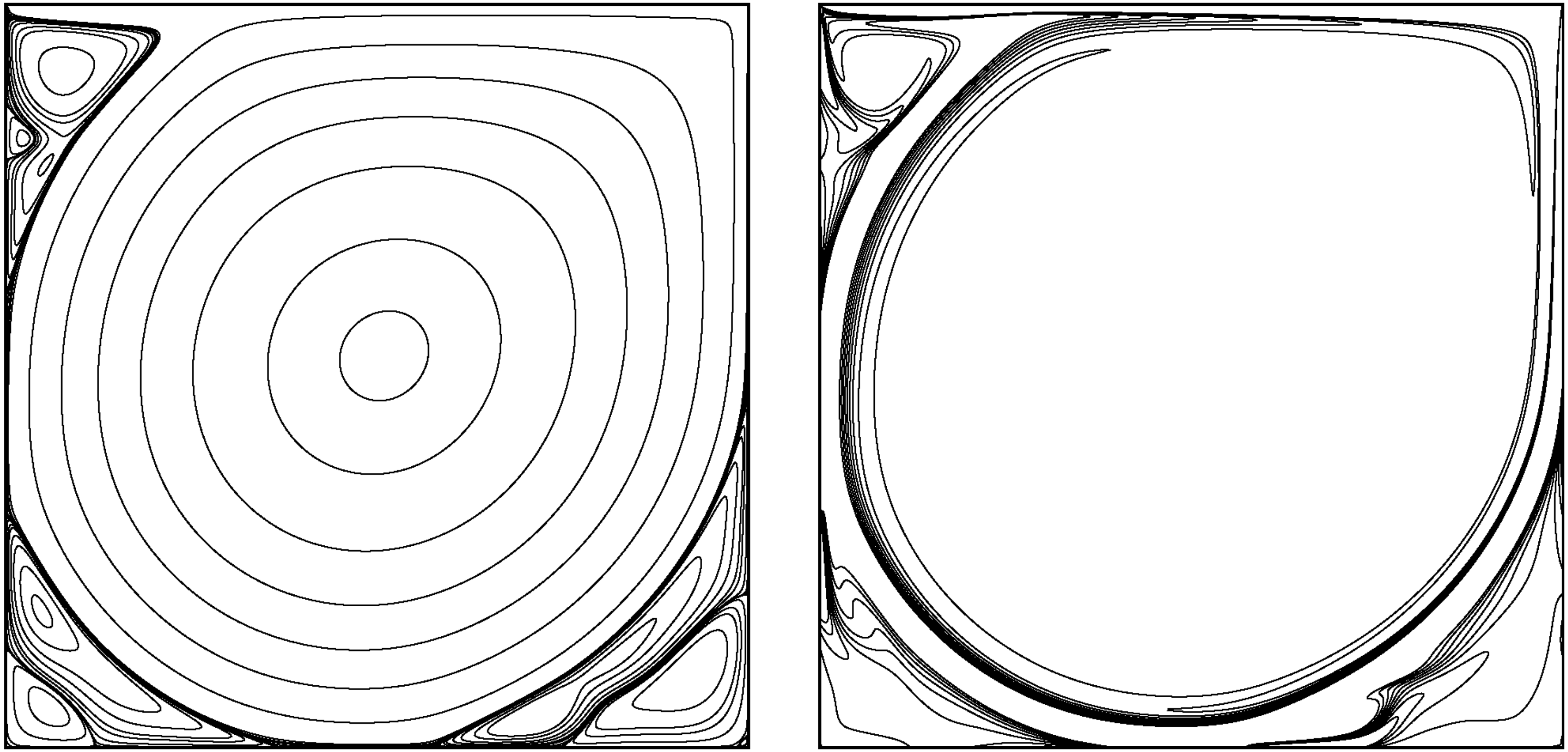

Figure 6. Streamfunction and vorticity contours for Re=20000

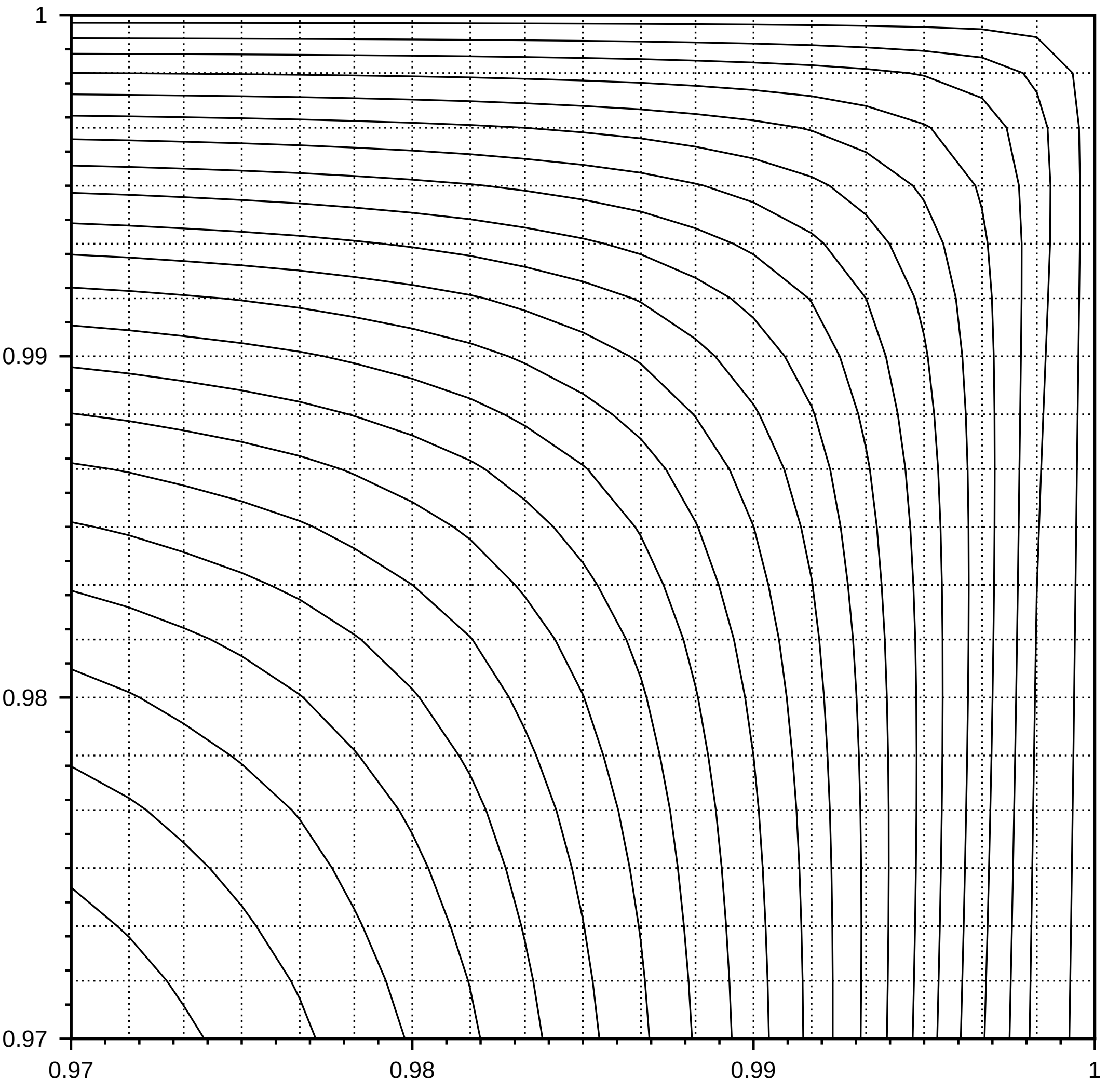

Figure 7. Streamfunction contours for Re=20000, enlarged view of top right corner.

| | |
|---|---|
| $\phi_x$ | $\frac{\phi_{i+1,j}-\phi_{i-1,j}}{2\Delta x}$ |
| $\phi_y$ | $\frac{\phi_{i,j+1}-\phi_{i,j-1}}{2\Delta y}$ |
| $\phi_{xx}$ | $\frac{\phi_{i+1,j}-2\phi_{i,j}+\phi_{i-1,j}}{\Delta x^2}$ |
| $\phi_{yy}$ | $\frac{\phi_{i,j+1}-2\phi_{i,j}+\phi_{i,j-1}}{\Delta y^2}$ |
| $\phi_{xy}$ | $\frac{\phi_{i+1,j+1}-\phi_{i-1,j+1}-\phi_{i+1,j-1}+\phi_{i-1,j-1}}{4\Delta x\Delta y}$ |
| $\phi_{xxy}$ | $\frac{\phi_{i+1,j+1}-2\phi_{i,j+1}+\phi_{i-1,j+1}-\phi_{i+1,j-1}+2\phi_{i,j-1}-\phi_{i-1,j-1}}{2\Delta x^2\Delta y}$ |
| $\phi_{xyy}$ | $\frac{\phi_{i+1,j+1}-2\phi_{i+1,j}+\phi_{i+1,j-1}-\phi_{i-1,j+1}-2\phi_{i-1,j}+\phi_{i-1,j-1}}{2\Delta x\Delta y^2}$ |
| $\phi_{xxyy}$ | $\frac{\phi_{i+1,j+1}-2\phi_{i,j+1}+\phi_{i-1,j+1}-2\phi_{i+1,j}+4\phi_{i,j}-2\phi_{i-1,j}+\phi_{i+1,j-1}-2\phi_{i,j-1}+\phi_{i-1,j-1}}{\Delta x^2\Delta y^2}$ |

Table 1. Standard Second Order Central Discretizations, $\mathcal{O}(\Delta x^2, \Delta y^2)$.

| | Re | 1000 | 2500 | 5000 | 7500 | 10000 | 12500 | 15000 | 17500 | 20000 |
|---|---|---|---|---|---|---|---|---|---|---|
| **Primary** | $\psi$ | -0.118938 | -0.121472 | -0.122216 | -0.122344 | -0.122306 | -0.122201 | -0.122060 | -0.121889 | -0.121694 |
| **Vortex** | $\omega$ | -2.067760 | -1.976132 | -1.940547 | -1.926478 | -1.918187 | -1.912307 | -1.907651 | -1.903659 | -1.900032 |
| | **(x , y)** | (0.5300 , 0.5650) | (0.5200 , 0.5433) | (0.5150 , 0.5350) | (0.5133 , 0.5317) | (0.5117 , 0.5300) | (0.5117 , 0.5283) | (0.5100 , 0.5283) | (0.5100 , 0.5283) | (0.5100 , 0.5267) |
| | $\psi$ | 0.17297E-02 | 0.26623E-02 | 0.30735E-02 | 0.32265E-02 | 0.31896E-02 | 0.30972E-02 | 0.30022E-02 | 0.29021E-02 | 0.28012E-02 |
| **BR1** | $\omega$ | 1.118222 | 1.955594 | 2.739071 | 3.243925 | 3.756425 | 4.357323 | 4.965304 | 5.568522 | 6.125275 |
| | **(x , y)** | (0.8633 , 0.1117) | (0.8333 , 0.0900) | (0.8050 , 0.0733) | (0.7900 , 0.0650) | (0.7750 , 0.0600) | (0.7600 , 0.0550) | (0.7450 , 0.0500) | (0.7333 , 0.0467) | (0.7200 , 0.0433) |
| | $\psi$ | 0.23345E-03 | 0.93093E-03 | 0.13758E-02 | 0.15337E-02 | 0.16135E-02 | 0.16568E-02 | 0.16663E-02 | 0.16450E-02 | 0.16083E-02 |
| **BL1** | $\omega$ | 0.354271 | 0.979831 | 1.514292 | 1.858636 | 2.182043 | 2.358513 | 2.490168 | 2.704118 | 2.885054 |
| | **(x , y)** | (0.0833 , 0.0783) | (0.0833 , 0.1117) | (0.0733 , 0.1367) | (0.0650 , 0.1517) | (0.0583 , 0.1633) | (0.0550 , 0.1683) | (0.0533 , 0.1717) | (0.0500 , 0.1767) | (0.0483 , 0.1817) |
| | $\psi$ | -0.49242E-07 | -0.12035E-06 | -0.14271E-05 | -0.32742E-04 | -0.14014E-03 | -0.25498E-03 | -0.34006E-03 | -0.40492E-03 | -0.46117E-03 |
| **BR2** | $\omega$ | -0.76887E-02 | -0.11935E-01 | -0.33647E-01 | -0.155445 | -0.309618 | -0.398536 | -0.456877 | -0.512862 | -0.554507 |
| | **(x , y)** | (0.9917 , 0.0067) | (0.9900 , 0.0083) | (0.9783 , 0.0183) | (0.9517 , 0.0417) | (0.9350 , 0.0683) | (0.9283 , 0.0817) | (0.9267 , 0.0883) | (0.9283 , 0.0967) | (0.9300 , 0.1050) |
| | $\psi$ | -0.65156E-08 | -0.26763E-07 | -0.65524E-07 | -0.20189E-06 | -0.11017E-05 | -0.64493E-05 | -0.22287E-04 | -0.48296E-04 | -0.78077E-04 |
| **BL2** | $\omega$ | -0.29579E-02 | -0.93547E-02 | -0.12183E-01 | -0.16587E-01 | -0.30153E-01 | -0.74310E-01 | -0.140685 | -0.200227 | -0.243259 |
| | **(x , y)** | (0.0050 , 0.0050) | (0.0067 , 0.0067) | (0.0083 , 0.0083) | (0.0117 , 0.0117) | (0.0167 , 0.0200) | (0.0267 , 0.0317) | (0.0383 , 0.0417) | (0.0500 , 0.0483) | (0.0583 , 0.0533) |
| | $\psi$ | - | 0.34284E-03 | 0.14459E-02 | 0.21311E-02 | 0.26237E-02 | 0.29949E-02 | 0.32865E-02 | 0.35245E-02 | 0.37247E-02 |
| **TL1** | $\omega$ | - | 1.344040 | 2.115544 | 2.236986 | 2.322310 | 2.375671 | 2.424136 | 2.463058 | 2.497615 |
| | **(x , y)** | - | (0.0433 , 0.8900) | (0.0633 , 0.9100) | (0.0667 , 0.9117) | (0.0700 , 0.9100) | (0.0733 , 0.9100) | (0.0767 , 0.9100) | (0.0783 , 0.9117) | (0.0800 , 0.9117) |
| | $\psi$ | - | - | 0.42053E-10 | 0.75831E-09 | 0.37451E-08 | 0.78212E-08 | 0.11742E-07 | 0.16951E-07 | 0.27181E-07 |
| **BR3** | $\omega$ | - | - | 0.66735E-03 | 0.70882E-03 | 0.24926E-02 | 0.37843E-02 | 0.33420E-02 | 0.42443E-02 | 0.49718E-02 |
| | **(x , y)** | - | - | (0.9983 , 0.0017) | (0.9967 , 0.0017) | (0.9967 , 0.0050) | (0.9950 , 0.0050) | (0.9950 , 0.0050) | (0.9950 , 0.0067) | (0.9933 , 0.0067) |
| | $\psi$ | - | - | - | - | - | 0.24668E-09 | 0.57300E-09 | 0.13069E-08 | 0.23598E-08 |
| **BL3** | $\omega$ | - | - | - | - | - | 0.50277E-03 | 0.10067E-02 | 0.21872E-02 | 0.18079E-02 |
| | **(x , y)** | - | - | - | - | - | (0.0017 , 0.0017) | (0.0017 , 0.0033) | (0.0033 , 0.0033) | (0.0033 , 0.0033) |
| | $\psi$ | - | - | - | - | - | -0.16253E-05 | -0.15848E-04 | -0.41183E-04 | -0.70762E-04 |
| **TL2** | $\omega$ | - | - | - | - | - | -0.204044 | -0.524342 | -0.737721 | -0.982968 |
| | **(x , y)** | - | - | - | - | - | (0.0067 , 0.8300) | (0.0150 , 0.8250) | (0.0200 , 0.8217) | (0.0250 , 0.8200) |

Table 2. Properties of primary and secondary vortices; streamfunction and vorticity values, (x,y) locations

| Re | present<br>$\psi_{min}$ ($\mathcal{O}\Delta x^4$) | Erturk *et. al.* [2]<br>$\psi_{min}$ ($\mathcal{O}\Delta x^4$) |
|---|---|---|
| 1000 | -0.118938 | -0.118939 |
| 2500 | -0.121472 | -0.121469 |
| 5000 | -0.122216 | -0.122213 |
| 7500 | -0.122344 | -0.122341 |
| 10000 | -0.122306 | -0.122313 |
| 12500 | -0.122201 | -0.122229 |
| 15000 | -0.122060 | -0.122124 |
| 17500 | -0.121889 | -0.122016 |
| 20000 | -0.121694 | -0.121901 |

Table 3. Minimum Streamfunction values at the Primary Vortex for various Reynolds numbers

| Re | present<br>$\omega$ ($\mathcal{O}\Delta x^4$) | Erturk *et. al.* [2]<br>$\omega$ ($\mathcal{O}\Delta x^4$) |
|---|---|---|
| 1000 | -2.067760 | -2.067579 |
| 2500 | -1.976132 | -1.976096 |
| 5000 | -1.940547 | -1.940451 |
| 7500 | -1.926478 | -1.926282 |
| 10000 | -1.918187 | -1.917919 |
| 12500 | -1.912307 | -1.912072 |
| 15000 | -1.907651 | -1.907602 |
| 17500 | -1.903659 | -1.903975 |
| 20000 | -1.900032 | -1.900891 |

Table 4. Vorticity values at the center of the Primary Vortex for various Reynolds numbers